# Women 1.5 Times More Likely to Leave STEM Pipeline After Calculus Compared to Men: Lack of Mathematical Confidence a Potential Culprit

Jess Ellis[1*], Bailey K. Fosdick[2], Chris Rasmussen[3]

[1]Department of Mathematics, Colorado State University, Fort Collins, Colorado, United States of America

[2]Department of Statistics, Colorado State University, Fort Collins, Colorado, United States of America

[3]Department of Mathematics and Statistics, San Diego State University, San Diego, California, United States of America

*Corresponding author:
E-mail: ellis@math.colostate.edu

# Abstract

The substantial gender gap in the science, technology, engineering, and mathematics (STEM) workforce can be traced back to the underrepresentation of women at various milestones in the career pathway. Calculus is a necessary step in this pathway and has been shown to often dissuade people from pursuing STEM fields. We examine the characteristics of students who begin college interested in STEM and either persist or switch out of the calculus sequence after taking Calculus I, and hence either continue to pursue a STEM major or are dissuaded from STEM disciplines. The data come from a unique, national survey focused on mainstream college calculus. Our analyses show that, while controlling for academic preparedness, career intentions, and instruction, the odds of a woman being dissuaded from continuing in calculus is 1.5 times greater than that for a man. Furthermore, women report they do not understand the course material well enough to continue significantly more often than men. When comparing women and men with above-average mathematical abilities and preparedness, we find women start and end the term with significantly lower mathematical confidence than men. This suggests a lack of mathematical confidence, rather than a lack of mathematically ability, may be responsible for the high departure rate of women. While it would be ideal to increase interest and participation of women in STEM at all stages of their careers, our findings indicate that simply increasing the retention of women starting in college calculus would almost double the number of women entering the STEM workforce.

# Introduction

Across the world there is tremendous need for more workers with degrees in science, technology, engineering, or mathematics (STEM). The U.S. President's Council of Advisors on Science and Technology (PCAST) report predicts over the next decade approximately one million more STEM graduates above and beyond the current graduation level will be needed in order to meet the demands of the U.S. workplace (*1*). The report also argues that simply increasing the retention of STEM majors by 10% would make considerable progress towards meeting this need.

In the United States and elsewhere, first-year college and university mathematics courses often function as a bottleneck, preventing large numbers of students from pursuing a STEM career (*2-4*). Introductory math courses, such as calculus, have repeatedly been linked to students' decisions to leave STEM majors (*5-7*). While calculus is not the only hurdle faced by potential U.S. STEM graduates, it is both one of the most challenging obstacles and a necessary first step on the way to a STEM career.

There has been a growing body of work investigating student persistence in STEM (*2-7*). A common perception is that students leave STEM majors because of poor academic ability and that calculus functions as a course that "weeds out" mathematically incapable students (*5, 8*). However, research suggests that switching from a STEM major to a non-STEM major is not an event, but a process based on a collection of curricular, instructional, and cultural issues (*9*). Seymour and her colleagues identified a number of these issues, including conceptual difficulties, poor instruction, inadequate preparation, and language barriers (*9, 10*). More recent work suggests that student demographics and socioeconomic status, secondary school preparation, student supports once in college, college pedagogy, and college grades are also important factors in STEM persistence (*11-14*).

In addition to this established "leaking STEM pipeline", women are underrepresented in STEM across all career stages. Although fourth-grade boys and girls report similar rates of interest in science, by twelfth-grade 34% of women and 48% of men report such an interest (*15*). By the time students enter college, 22% of women intend to study a STEM field compared to 34% of men (*16*). An estimated 40-60% of students who begin a STEM degree actually complete one, and of those only 29% are completed by women (*16, 17*). Combined, these decreases in women's participation in STEM lead to women making up only 25% of STEM workforce (*18*) (Fig. 1). In looking specifically within academia, these patterns persist, and although more women are entering academic positions than before, women continue to be an underrepresented minority in many STEM fields (*19, 20*). Studies indicate that while there exists no bias against women in hiring for tenure track positions (*21*), women are not afforded the same opportunities, such as elite post doc positions, that men are that help them be attractive for top academic positions (*22*).

As the U.S. faces a STEM graduate deficit, it is critical we understand why women and men are not completing STEM degrees at comparable rates and why both genders are not persisting with STEM degrees. In this study, we examine the role of Calculus I in STEM persistence for all students, focusing specifically on the gender gap. If a student elects not to take Calculus II, he or she is effectively choosing to exit the STEM pipeline. Thus, intentions to continue studying calculus after Calculus I may serve as a proxy for continuing to study STEM.

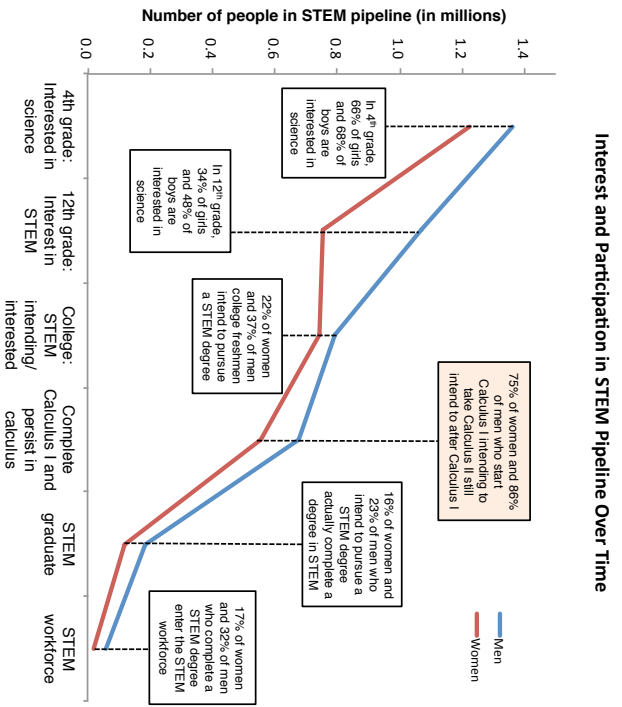
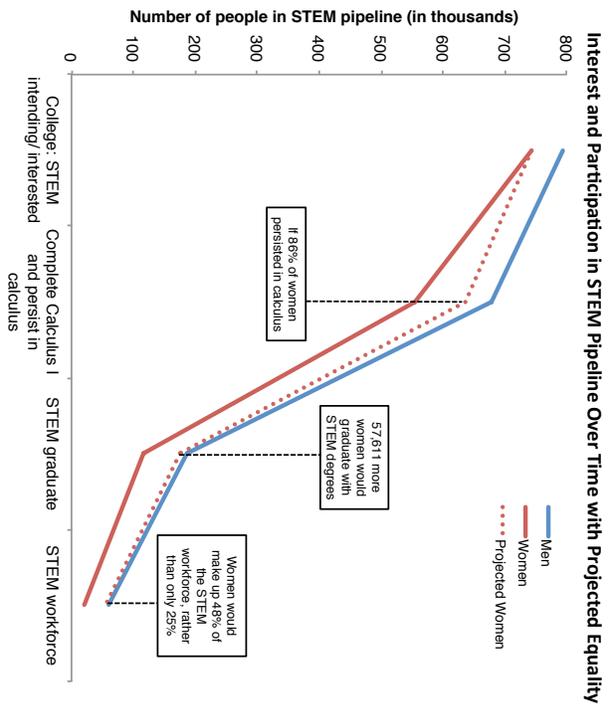

Fig 1. Actual participation of women in STEM (A) and projected participation by women in STEM if they decreased participation at the same rates as men after entering college (B). [n=3111]

# Current Study

The data used for this study comes from a unique, large-scale and in-depth national survey of Calculus I conducted under the auspices of the Mathematical Association of America. Colleges and universities were selected to participate using a stratified random sample of two- and four-year undergraduate colleges and universities during the 2010 Fall term. The San Diego State University Institutional Review Board (IRB) approved the study. The protocol number is 496064. Participant responses were de-identified prior to analysis.

The surveys were constructed based on a literature review of potential factors related to student success in calculus and feedback from experts on the projects' advisory board. Administration of the surveys were restricted to what is known as "mainstream" calculus, the calculus course designed to prepare students for studying engineering or the physical sciences. Until now, there has been very little large-scale data collected on who elects to study Calculus I or on the effect of this course on student persistence in STEM.

Students were surveyed at the beginning and end of the Calculus I term and asked if they intended to take Calculus II. One year later, students were asked if they had taken or enrolled in Calculus II. Based on students' responses, we identified students who initially intended to take Calculus II and noted whether this intent was maintained or not after Calculus I (see Table S1 for more information). Those who maintained their initial intention to take more calculus are referred to as *Persisters* and those who reported lower intentions of taking Calculus II at the end of term compared to the beginning of the term are referred to as *Switchers*.

In this study we examine the characteristics of students who enroll in Calculus I and either persist or switch out of the mainstream calculus sequence, and hence either remain or leave the STEM pipeline, attending specifically to gender. We perform a statistical analysis of student change in their intention to take Calculus II by gender, comparing Switchers to Persisters, while controlling for students' preparedness for Calculus I, intended career goals, perception of instruction, and institutional environment.

To measure preparedness, we use student reported previous calculus experience and standardized math test score (ACT and SAT). Career goals are characterized by students' reported career aspirations. Students intending to pursue a career in science, technology, or math are grouped together and labeled STM. We consider students pursuing medical professions, non-STEM fields (e.g. business, law, education), and those who are undecided to be STEM-interested as these students indicated they were originally planning to take Calculus II at the beginning of the term, and thus must have been initially open to pursing a degree that required more mathematics (see Table S2 for more detailed information). The STEM-interested students could be considering a STEM field as a second degree or interdisciplinary studies involving STEM – fields which are witnessing much greater demands in industry than specialized science fields *(23)*.

Student perception of instruction was characterized by aggregate variables *Instructor Quality* and *Student-Centered Practices*, ranging from 1-6, based on student reports of sixteen instructional practices and behaviors (see Tables S3 and S4 for detailed information on the derivation of these variables). Instructor Quality characterizes the level of conventional quality teaching, including availability outside of office hours, listening to questions, and encouraging students mathematically. Low values on this scale indicate low perceived instructional quality, and high values correspond to high instructional quality. *Student-Centered Practices* characterizes the frequency of classroom practices such as whole-class discussion, students giving presentations, and group work. Low values coincide with traditional, instructor-centered instructional practices, and high values correspond to more innovative, student-centered teaching. Since there are many unknown and unmeasured characteristics of an institutional environment that likely contribute to a student's career decision, we expect dependence among switcher propensities at the same institution. For this reason, we also included institution in the analysis.

# Results

## Identifying who is Switching out of Calculus

Using a logistic mixed effects regression model, we analyzed the association between switcher propensity and gender, controlling for student preparedness, career intentions, instruction, and institution. There were 2,266 students for which we had complete data and of these 17.8% were identified as Switchers.

Three of the controlling variables were found to be significant when predicting persistence: standardized math score, career intentions and Instructor Quality. As shown in Figure 2, higher standardized test scores correspond to an increased likelihood of persisting, as does intending to be an Engineer (compared to a STM field) and higher levels of Instructor Quality. Compared to students pursuing a STM field, students pursuing a medical field, non-STEM field, or are undecided are more likely to switch out of calculus. Interestingly, neither previous calculus experience nor Student-Centered Practices are significantly associated with switching propensity (since the credible intervals for the odds ratios contain one).

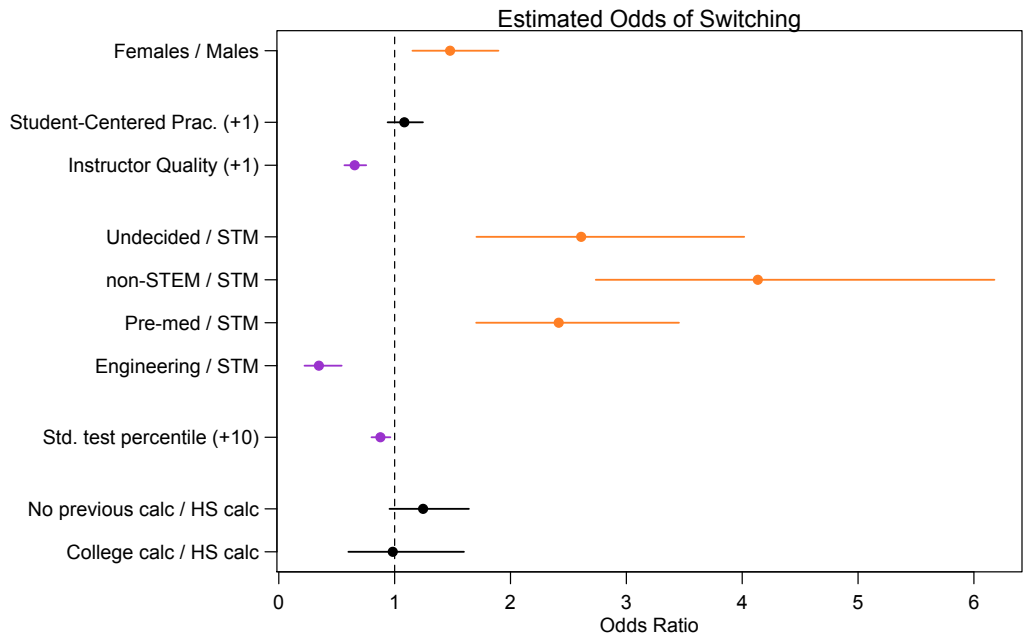

**Fig 2. Odds ratios of switching for student attributes.** The circle represents the odds ratio estimate and the bars represent the 95% credible interval. The continuous variables noted with (+*x*) on the left compare a student who reported *x*-points higher than another student. Labels of the form *A / B* correspond to the ratio of the odds of switching for a student of type *A* to the odds of switching for a student of type *B*. Variables associated with decreased likeliness and increased likeliness of switching are highlighted in purple and orange, respectively. [n=2266]

Even after controlling for student preparedness, career intentions, and instruction, gender is significantly related to persistence. Specifically, a female student's odds of switching are approximately 1.5 times that of a comparable male student of the same preparedness, career goals, and reports of instruction (95% CI: 1.15-1.90) (see Table S9 for more detail). To understand what this means practically, consider two hypothetical students, one STEM-intending and one STEM-interested: Student A earned an average standardized math score, took high school calculus, is pursuing STM, and reports average levels of Student-Centered Practices and lower than average levels of Instructor Quality. If student A is a man, he has an 11.8% probability of switching out of his calculus, whereas if student A is a woman, this probability increases to 16.5%. Student B also earned an average standardized math score, did not take high school calculus, is pursuing a non-STEM career, and reports average Instructor Quality and Student-Centered Practices. If student B is male, he has a 31.1% probability of switching out of his calculus, however if instead, student B is female, this probability increases to 40.0%. These results show that Calculus I is a critical "leak" in the STEM pipeline, especially for women.

## Examining Students' Reasons for Leaving Calculus

We now consider the question of why. On the end of term survey, students who did not intend to take Calculus II were given a list of potential reasons and were asked to select all that resonated with them. In Table 1, we report statistics on the reasons Switchers gave for not persisting in calculus at the end of Calculus I. These students represent a small sample of Switchers in our study, but their opinions provide insight into the potential beliefs of other students who did not share their reasons.

**Table 1. Switchers' reasons for not intending to take Calculus II.**

| | STEM-Intending | | STEM-Interested | |
|---|---|---|---|---|
| Reason for not intending to take Calc. II | Men (37) | Women (48) | Men (86) | Women (158) |
| I changed my major and now do not need to take Calculus II | 70% | 65% | 33% | 32% |
| To do well in Calculus II, I would need to spend more time and effort than I can afford | 41% | 35% | 38% | 37% |
| My experience in Calculus I made me decide not to take Calculus II | 32% | 38% | 42% | 45% |
| I have too many other courses I need to complete | 27% | 25% | 50% | 50% |
| I do not believe I understand the ideas of Calculus I well enough to take Calculus II | 14% | 35% | 20% | 32% |
| My grade in Calculus I was not good enough for me to continue to Calculus II | 16% | 19% | 15% | 15% |

Note: Gender differences that are statistically significant at the 0.10 level based on Fisher's exact test are highlighted in grey. The corresponding p-values for STEM-intending and STEM-interested students are 0.026 and 0.051, respectively. [n=329]

Common reasons selected by all students were a change in major, too many other courses to complete, their experience in Calculus I and the perception that Calculus II would require excessive time and effort. The proportions of students who cited each reason were comparable across men and women, except for one: "I do not believe I understand the ideas of Calculus I well enough to take Calculus II." Among STEM-intending students, 35% of women reported this as a reason while only 14% of men acknowledged it (p = 0.026). Among STEM-interested students, 32% of women reported this as a reason compared to only 20% of men (p = 0.051). Thus, women Switchers are citing a lack of understanding of the material in Calculus I as a reason for not continuing their studies significantly more often than men.

Previous research suggests that this is not because women do not actually understand the material as well as men; on the contrary, a meta-analysis of gender differences in mathematics found no differences in ability (*24*) and a study specifically looking at gender differences in Calculus I found that women outperform men (*25*). These gender differences are disconcerting as they suggest that perception of one's ability plays a role in women's decisions to stop taking calculus but not as much for men.

## Investigating Confidence as a Source for Gender Disparities in Switching

It is well documented that confidence plays a significant role in one's success (*26*), and that men and women have different levels of confidence in their mathematical ability (*27, 28*). This begs the question of whether calculus is weeding out students based on capability or a lack of confidence in their mathematical capability.

To explore this question, we compare the change in student reported mathematical confidence among mathematically-capable students grouped by gender and persistence. We operationalize mathematically-capable as those students with standardized math scores at or above the national 85[th] percentile. Figure 3 shows that all mathematically-capable students lose mathematical confidence over the course of Calculus I. Switchers experience a greater decrease in confidence than Persisters, and women start at a lower confidence and therefore end at a lower confidence, while experiencing the same decrease as men.

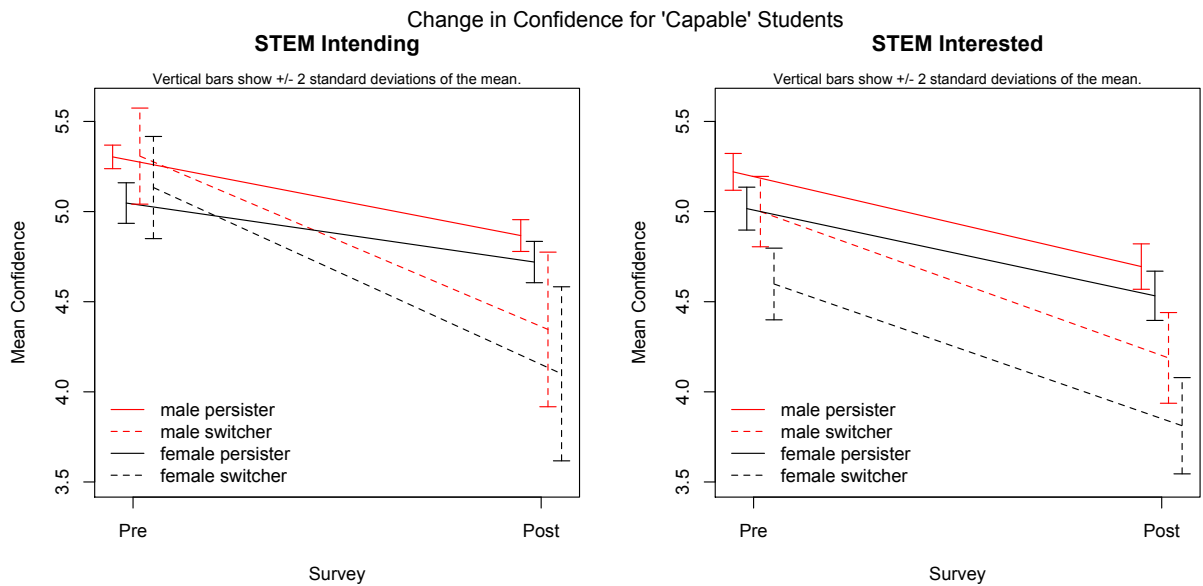

**Fig 3**. **Change in student mathematical confidence at the beginning of the Calculus I semester (pre-survey) and at the end of the semester (post-survey) separated by career intentions, gender, and persistence status.** [n=1524]

## Discussion

Calculus I is an established milestone in the STEM trajectory, and we have shown here that it is contributing significantly to the STEM "gender filter" (*29*). What can we do with this information? Our work points to women's mathematical confidence as a major factor in their decision not to persist in calculus, and therefore STEM. While men and women lose confidence at similar rates during Calculus I, they come into college calculus with different levels of mathematical confidence. Returning to Fig. 1, we see the transition from 4[th] grade to 12[th] grade as a critical decrease in women's interest in science. There is a clear need to target efforts at this age group. However, strong gains can be made even if women continue to enter college with lower levels of scientific interest and mathematical confidence, as shown by the dotted line. If women persisted in STEM at the same rate as men starting in Calculus I, women would make up 48% of the STEM workforce rather than the current 25%. Certainly it is preferable to increase girl's and women's interest in STEM at all life stages, but this projection indicates that only targeting efforts at college calculus and beyond would almost double the number of women entering the STEM workforce. This would increase the incoming STEM workforce by 45%, and go a long way to meet the needs articulated in the PCAST report (*1*).

## Acknowledgments

This work is part of the *Characteristics of Successful Programs in College Calculus* project and would not have been possible without the many contributors to this work, especially David Bressoud. The authors are grateful to feedback from Sandra Laursen, Jennifer Hoeting, Tyler McCormick, and Rachel Heath.

# Supplementary Information: Methods

## Data Preparation

*Switcher coding*

Students were coded as a Switcher, Persister, or neither based on their responses to four questions. On the beginning of term survey, students were asked if they intended to take Calculus II, with options "yes", "no" or "I don't know." On the end of term survey, students were asked if, at the beginning of the term, they intended to take Calculus II, with options "yes", "no", or "I don't remember" (referred to as "End of term; reflect"). They were also asked if they currently intended to take Calculus II, with options "yes", "no", or "I don't know." On the follow up survey on year later, students were asked if they had already taken or were currently enrolled in Calculus II, with options "yes" or "no."

Because not all students answered all surveys, there were multiple ways that we identified students as Switchers, Persisters, or neither. Switchers were identified as any student who gave sufficient evidence of decreasing their intentions to take Calculus II. There were 11 unique ways that students could be identified as a Switcher, as shown in Table S1. In this table, student responses to each of the four questions are filled in as "Y" for "yes", "N" for "no", "M" for "I don't know" or "I don't remember", "NA" for not answered, and blank for any option. Thus, in Switcher group 1, students answered "yes" to the beginning of term question, answered anything to the end of term survey questions, and answered "no" a year later. Thus, this group of students entered Calculus I intended to take Calculus II and a year later had not taken Calculus II. Students in Switcher group 8 were initially unsure whether they would take Calculus II, marking "I don't know". However, this uncertainty suggests they were at least interested or open to taking more calculus. By the end of the term they said that they did not intend to take Calculus II, and thus they decreased their STEM interest and/or intention.

Students whose responses are not captured in the table below were determined to not be initially interested or intending to take Calculus II.

**Table S1.** Switcher coding dictionary

| Switchers | Number | 1 – Beginning of term | 2 – End of term; reflect | 3 – End of term | 4 - Follow Up |
|---|---|---|---|---|---|
| 1 | 160 | Y | | | N |
| 2 | 118 | M | | | N |
| 3 | 15 | NA | Y | | N |
| 4 | 3 | NA | M | | N |
| 5 | 38 | Y | Y | M | NA |
| 6 | 123 | Y | | N | NA |
| 7 | 17 | M | Y | M | NA |
| 8 | 152 | M | | N | NA |
| 9 | 34 | NA | Y | M | NA |
| 10 | 78 | NA | Y | N | NA |
| 11 | 65 | NA | M | N | NA |

| Persisters | Number | 1 – Beginning of term | 2 – End of term; reflect | 3 – End of term | 4 - Follow Up |
|---|---|---|---|---|---|
| 12 | 586 | Y | | | Y |
| 13 | 63 | M | | | Y |
| 14 | 67 | NA | Y | | Y |
| 15 | 2 | NA | M | | Y |
| 16 | 1543 | Y | | Y | NA |
| 17 | 35 | Y | M | M | NA |
| 18 | 5 | Y | N | M | NA |
| 19 | 1 | Y | NA | M | NA |
| 20 | 193 | M | | Y | NA |
| 21 | 64 | M | M | M | NA |
| 22 | 22 | M | N | M | NA |
| 23 | 3 | M | NA | M | NA |
| 24 | 1325 | NA | Y | Y | NA |
| 25 | 53 | NA | M | Y | NA |
| 26 | 103 | NA | M | M | NA |

*Career choice grouping*

On the beginning of term survey students were asked to indicate their intended career choice, choosing from one of 16 options, shown in Table S2. These intended careers were grouped together into five groups. The first group is comprised of traditional STEM degrees, excluding Engineering, the second group. We chose to exclude engineering because there were a disproportionate number of engineering students compared to the other STEM fields. The third group of career intentions is made up of medical and other health professionals. The fourth group is made up of traditionally non-STEM careers, including non-STEM education, social scientists, business, law, humanities, and other non-science related career. The final category is Undecided.

Rather than restrict our analysis to students who indicated that they are intending to pursue a career in STEM, we included all students who indicated that they were at least open to taking more calculus, and thus were either STEM-intending or STEM-interested.

**Table S2.** Career choice grouping

| Group | Coding | Original indicated career choice on beginning of term survey |
|---|---|---|
| STM – traditional STEM fields, excluding engineers | 1 | Life scientist (e.g. biologist, medical researcher) |
| | 1 | Earth/Environment scientist (e.g. geologist, meteorologist) |
| | 1 | Physical scientist (e.g. chemist, physicist, astronomer) |
| | 1 | Computer scientist, IT |
| | 1 | Mathematician |
| | 1 | Science/Math teacher |
| Engineering | 2 | Engineer |
| Pre-med | 3 | Medical professional (e.g. doctor, dentist, vet) |
| | 3 | Other health professional (e.g. nurse, medical technician) |
| Non-STEM | 4 | Other teacher |
| | 4 | Social scientist (e.g. psychologist, sociologist) |
| | 4 | Business administration (finance, accounting, management) |
| | 4 | Lawyer |
| | 4 | English/Language/Arts specialist |
| | 4 | Other non-science related career |

*Reports of instruction*

On the end of term survey, students were asked two multi-part questions related to instructional practices. Both questions were on a 6-point scale. For the first set of questions, 1 indicated strongly disagree and 6 indicated strongly agree. For the second set of questions, 1 indicated not at all, and 6 indicated very often.

**18. My Calculus instructor:**

|  | Strongly Disagree | Disagree | Slightly Disagree | Slightly Agree | Agree | Strongly Agree |
|---|---|---|---|---|---|---|
| asked questions to determine if I understood what was being discussed. | ○ | ○ | ○ | ○ | ○ | ○ |
| listened carefully to my questions and comments. | ○ | ○ | ○ | ○ | ○ | ○ |
| discussed applications of calculus. | ○ | ○ | ○ | ○ | ○ | ○ |
| allowed time for me to understand difficult ideas. | ○ | ○ | ○ | ○ | ○ | ○ |
| helped me become a better problem solver. | ○ | ○ | ○ | ○ | ○ | ○ |
| provided explanations that were understandable. | ○ | ○ | ○ | ○ | ○ | ○ |
| was available to make appointments outside of office hours, if needed. | ○ | ○ | ○ | ○ | ○ | ○ |
| discouraged me from wanting to continue taking Calculus. | ○ | ○ | ○ | ○ | ○ | ○ |

(a)

**19. During class time, how frequently did your instructor:**

|  | Not at all 1 | 2 | 3 | 4 | 5 | Very often 6 |
|---|---|---|---|---|---|---|
| show how to work specific problems? | ○ | ○ | ○ | ○ | ○ | ○ |
| have students work with one another? | ○ | ○ | ○ | ○ | ○ | ○ |
| hold a whole-class discussion? | ○ | ○ | ○ | ○ | ○ | ○ |
| have students give presentations? | ○ | ○ | ○ | ○ | ○ | ○ |
| have students work individually on problems or tasks? | ○ | ○ | ○ | ○ | ○ | ○ |
| lecture? | ○ | ○ | ○ | ○ | ○ | ○ |
| ask questions? | ○ | ○ | ○ | ○ | ○ | ○ |
| ask students to explain their thinking? | ○ | ○ | ○ | ○ | ○ | ○ |

(b)

**Figure S1**: Instructor course practices related to (a) "instructor quality" and (b) "student-centered practices".

Informed by a factor analysis of these 16 practices, the eight practices from question 18 (Figure S1) were averaged to create a new variable called *Instructor Quality,* with the last prompt (my instructor discouraged me from wanting to continue taking Calculus) reverse coded. Thus, for this new variable a 1 indicates low levels of Instructor Quality and 6 indicated high levels. The eight practices from question 19 (Figure S2) were averaged to create a new variable called *Student-Centered Practices* with two prompts (show how to work specific problems and lecture) weighted .5 based on the factor analysis. Again, for this new variable a 1 indicates low levels of Student-Centered Practices and 6 indicated high levels. Tables S3 and S4 show the loadings from the factor analysis, and the weights of the averages in creating the aggregate variables. We chose to use averages instead of the PCA loadings as weights so that the variables were more easily interpreted. For both variables, if students did not answer some of the questions we took the average of the other variables.

**Table S3.** PCA loadings and aggregate variable weights for *Instructor Quality*

|  | PCA loading | Aggregate variable weights |
|---|---|---|
| 18. My calculus instructor: (1 – strongly disagree; 6 – strongly agree) | | |
| Asked questions to determine if I understood what was being discussed | 0.37 | 0.125 |
| Listened carefully to my questions and comments | 0.36 | 0.125 |
| Discussed applications of calculus | 0.30 | 0.125 |
| Allowed time for me to understand difficult ideas | 0.41 | 0.125 |
| Helped me become a better problem solver | 0.41 | 0.125 |
| Provided explanations that were understandable | 0.41 | 0.125 |
| Was available to make appointments outside of office hours, if needed. | 0.23 | 0.125 |
| Discouraged me from wanting to continue taking Calculus. | -0.29 | -0.125 |

**Table S4.** PCA loadings and aggregate variable weights for *Student-Centered Practices*

|  | PCA loading | Aggregate variable weights |
|---|---|---|
| 19. During class time, how frequently did your instructor: (1 – not at all; 6 – very often) | | |
| Show how to work specific problems? | 0.09 | 0.0625 |
| Have students work with one another? | 0.50 | 0.125 |
| Hold a whole-class discussion? | 0.51 | 0.125 |
| Have students give presentations? | 0.29 | 0.125 |
| Have students work individually on problems or tasks? | 0.37 | 0.125 |
| Lecture? | 0.05 | 0.0625 |
| Ask questions? | 0.25 | 0.125 |
| Ask students to explain their thinking? | 0.46 | 0.125 |

*Mathematical preparation*

To measure students' mathematical preparation, we use their previous calculus experience and their math standardizes test scores. For previous calculus experience, we group previous experience in calculus into three bins: high school (non-AP, AP AB, or AP BC), college, and none. For standardized test scores, students were asked to report their SAT math test score and/or their ACT math test score. Using the college board website and the ACT website[1] reports of percentiles, we converted these scores to an aggregate "Standardized math test score" code. For students who reported both, we used the average of their percentiles. See Table S6 for a summary of these variables.

---

[1] https://secure-media.collegeboard.org/digitalServices/pdf/sat/sat-percentile-ranks-crit-reading-math-writing-2014.pdf; http://www.actstudent.org/scores/norms1.html

**Descriptive analysis of switchers**

In the Tables S5-S8 below we summarize the relationships between student covariates, switcher code and gender.

**Table S5.** Descriptive table of career choice, by switcher code and gender

| Career Choice | Switcher | | Total Male | Total Female |
| --- | --- | --- | --- | --- |
| | *Male (N=166)* | *Female (N=238)* | *1236* | *1030* |
| STM | 10.6% | 16.1% | 263 | 223 |
| Engineering | 3.3% | 6.4% | 538 | 249 |
| Pre-med | 21.6% | 33.3% | 199 | 318 |
| Non-STEM | 36.8% | 38.1% | 136 | 126 |
| Undecided | 26.3% | 28.1% | 99 | 114 |

**Table S6.** Descriptive table of preparation, by switcher code and gender

| | Persister | | Switcher | |
| --- | --- | --- | --- | --- |
| **Previous Calculus** | *Male (N=1070)* | *Female (N=792)* | *Male (N=166)* | *Female (N=238)* |
| High School | 644 | 547 | 96 | 152 |
| College | 82 | 53 | 17 | 9 |
| None | 344 | 192 | 53 | 77 |
| **Standardized math test score *(Percentile)*** | | | | |
| 100-90% | 591 | 374 | 84 | 112 |
| 90-80% | 241 | 200 | 41 | 54 |
| 80-70% | 135 | 107 | 21 | 39 |
| 70-60% | 51 | 60 | 8 | 13 |
| 60-50% | 20 | 21 | 6 | 8 |
| 50-40% | 17 | 24 | 4 | 7 |
| 40-30% | 5 | 4 | 1 | 2 |
| 30-20% | 6 | 1 | 0 | 1 |
| 20-10% | 1 | 1 | 1 | 1 |
| 10-0% | 3 | 0 | 0 | 1 |
| Ave. (std. dev.) | 86.79 (13.56) | 85.14 (13.58) | 85.66 (14.11) | 83.81 (15.45) |

**Table S7.** Descriptive table of perceptions of instruction, by switcher code and gender

|  | Persister | | Switcher | |
|---|---|---|---|---|
| **Instructor Quality** | Male (N=1070) | Female (N=792) | Male (N=166) | Female (N=238) |
| 1-1.5 | 4 | 2 | 1 | 1 |
| 1.5-2.5 | 22 | 17 | 6 | 18 |
| 2.5-3.5 | 63 | 50 | 19 | 22 |
| 3.5-4.5 | 259 | 189 | 45 | 69 |
| 4.5-5.5 | 550 | 387 | 77 | 101 |
| 5.5-6 | 172 | 147 | 18 | 27 |
| Ave. (std. dev.) | 4.66 (0.87) | 4.69 (0.88) | 4.42 (0.96) | 4.32 (1.06) |

| **Student-Centered Practices** | Male (N=1070) | Female (N=792) | Male (N=166) | Female (N=238) |
|---|---|---|---|---|
| 1-1.5 | 30 | 30 | 7 | 7 |
| 1.5-2.5 | 217 | 179 | 23 | 54 |
| 2.5-3.5 | 344 | 252 | 67 | 77 |
| 3.5-4.5 | 331 | 228 | 46 | 62 |
| 4.5-5.5 | 132 | 90 | 19 | 29 |
| 5.5-6 | 16 | 13 | 4 | 9 |
| Ave. (std. dev.) | 3.26 (1.03) | 3.17 (1.06) | 3.26 (1.03) | 3.26 (1.11) |

**Table S8.** Pre- and Post-term reports of confidence, by switch code and gender

|  | Persister | | Switcher | |
|---|---|---|---|---|
|  | Male (N=1067) | Female (N=790) | Male (N=165) | Female (N=237) |
| Pre-term | 5.119 | 4.887 | 4.914 | 4.641 |
| Post-term | 4.709 | 4.474 | 4.251 | 3.825 |
| Difference | 0.41 | 0.413 | 0.663 | 0.816 |

**Statistical Analysis of Switchers**

A logistic regression model was used to quantify the association between various student characteristics and the propensity for students to switch out of the mainstream calculus sequence. Student standardized test score, previous calculus experience, career goals, course teaching perceptions, and gender were treated as fixed effects and institution was treated as a random effect. Career goals, previous calculus experience, gender and institution are categorical variables, while all standardized test score, good and Student-Centered Practices are continuous. Parameter estimates were obtained using Bayesian methods, where prior distributions were specified for all parameters and inference was based on the posterior distribution of the parameters given the data. An approximation to the posterior distribution was obtained using standard Markov chain Monte Carlo techniques.

Weakly-informative prior distributions were specified for the parameters. Specifically the prior distributions on the regression coefficients were independent mean-zero normal distributions with a common variance of 100. The prior distributions on the random effects were independent normal distributions with mean zero and a common variance parameter. The hyper-prior distribution on the variance parameter was a diffuse inverse-gamma distribution.

The chain was run for five million iterations, with an additional ten thousand burn-in iterations. Samples from the posterior were collected every $500^{th}$ iteration to reduce dependence in the samples. This resulted in 10,000 samples from the posterior distribution of the parameters given the data, where the effective sample sizes for all parameters was greater than 578. The posterior mean parameter estimates resulting from the MCMC procedure were compared with the maximum likelihood estimates given by the glmer function in the lme4 package in R (*30*), based on an adaptive Gauss-Hermite quadrature procedure with ten points per axis, and those from GLIMMIX in SAS. The estimates from all three estimation procedures were similar. The MCMC procedure was initialized with the estimates from glmer.

Summaries of the posterior distribution of the parameters are given in Figure 2 in the manuscript and Table S9. The point estimates are the means of the posterior distribution and the 95% credible intervals were created from the $2.5^{th}$ and $97.5^{th}$ quantiles. The posterior mean estimate of the intercept (on the logit scale) and 95% credible interval are -2.43 and (-2.82, -2.07). Also, the estimate and credible interval for the variance of the institution random effects are 0.49 and (0.25, 0.86).

**Table S9.** Logistic mixed effects model summary. Odds ratios for switching for categorical variables are presented relative to the reference category noted next to the characteristic (e.g. odds ratio for switching for college calculus is for college calculus compared to high school calculus). The odds ratio for switching for standardized test score compares a student with a test score 10 percentiles higher than another comparable student. Instructor Quality and Student-Centered Practices odds ratios compare perceived instruction for a student rating the course 1 unit higher than another student. Effects with odds ratio credible intervals (CI) that do not contain one are considered to be significant predictors of persistence.

|  | odds ratio | 95% CI |
|---|---|---|
| Previous calculus: compared to HS calculus | | |
| College calculus | 0.984 | (0.601, 1.598) |
| None | 1.246 | (0.956, 1.641) |
| | | |
| Standardized test score: | | |
| Percentile (10 pt increase) | 0.877 | (0.800, 0.963) |
| | | |
| Career choice: compared to STM | | |
| Engineering | 0.346 | (0.222, 0.543) |
| Pre-med | 2.416 | (1.704, 3.454) |
| non-STEM | 4.135 | (2.736, 6.177) |
| Undecided | 2.610 | (1.706, 4.017) |
| | | |
| Reports of instruction: | | |
| Instructor Quality (1 pt increase) | 0.655 | (0.566, 0.754) |
| Student-Centered Practices (1 pt increase) | 1.084 | (0.940, 1.244) |
| | | |
| Gender: compared to men | | |
| Women | 1.478 | (1.154, 1.896) |

**Limitations**

Due to the large-scale and complex nature of the study, there are a number of limitations of our work. First, we have no measure of how the students performed in their calculus courses. Grades were collected for only a very small percentage of students and the study involved no other measures of student calculus performance. The second limitation is in potential non-response bias. The survey was sent to a stratified random sample of institutions. Within those institutions, it was sent to all Calculus I instructors, who were instructed to send it to all of their Calculus I students at the beginning and at the end of the Calculus I term. There were multiple opportunities within this structure for both instructors and students to choose to opt out of participating in the study. For example, if a student strongly disliked Calculus I and dropped out of the class before the end of the term, he or she may have not filled out the end-of-term survey. In this case, although the student is a Switcher, he or she would not be included in our study. We project that these non-responses would lead to a higher percentage of Switchers than reported here. Finally, since we do not have data for the students who did not respond to the surveys, we cannot compare the students in our sample to the general population.